# A Flow-based Distributed Trading Mechanism in Regional Electricity Market with Energy Hub

Lu Wang [1, 3], Mokhtar Bozorg [2], Mohammad Rayati [2], and Rachid Cherkaoui [1]

[1] EPFL, 1015 Lausanne, Switzerland,
{lu.wang,rachid.cherkaoui}@epfl.ch

[2] HES-SO University of Applied Sciences of Western Switzerland, Route de Cheseaux 1, 1401 Yverdon-les-Bains, Switzerland
{ mokhtar.bozorg, Mohammad.rayati }@heig-vd.ch

[3] Southeast University, Nanjing, Jiangsu Province, 210096 China

*Abstract*—The concept of Energy Hub (EH) has been emerged to accommodate renewable energy sources in a multi-energy system to deploy the synergies between electricity and other energy sources. However, the market mechanisms for integration of the EHs into the energy markets are not sufficiently elaborated. This paper proposes a flow-based two-level distributed trading mechanism in the regional electricity market with EH. At the lower level, the regional system operator coordinates the regional grids transactions in two markets, the local energy market with EH and the wholesale market of the upstream grid. Every nodal agent as an independent stakeholder leverages price discrepancy to cross arbitrage from different markets. At the upper level, the EH is a third player intending to maximize profit from trading in the regional electricity market and gas market. The regional electricity market clearing problem is formulated as a mathematical program with equilibrium constraints, for which we develop an ADMM-based distributed algorithm to obtain the equilibrium solution. The DC power flow is decomposed into optimization problems for the regional system operator and agents at different nodes, which can be solved in a distributed manner to achieve global optimality without violating the privacy of players. Case studies based on a realistic regional grid verify the effectiveness of the proposed algorithm and show that the mechanism is effective in decomposing power flow and incrasing energy efficiency.

## I. Introduction

The sustainable development of human society is restricted by climate change and energy shortage [1]. Many countries are encouraging the utilization of renewable energy sources while promoting multi-energy coupling. From technical perspective, the concept of Energy Hub (EH) has been emerged to accommodate renewable energy sources in a multi-energy system to deploy the synergies between electric power and other energy sources [2]. From a market perspective, the establishment of EHs and distributed energy sources requires an efficient trading mechanism in regional electricity markets.

To adapt to the development of distributed energy resources and multi-energy coupling, Li et al. [3] proposed a Lyapunov-based energy management method to achieve economic operation among EHs. An electricity-heat retail market framework was established in [4] to realize the optimal energy allocation of energy station devices. A comprehensive optimal bidding strategy for an energy hub was modeled in [5], which enabled the energy hub to benefit from day-ahead and real-time market. In short, the previous studies, which profound techniques for performing trading mechanisms in a centralized manner, do not consider an important point that every stakeholder has a sense of privacy information protection to reveal their true cost or value functions [6]. To address this point, we need a distributed algorithm for implementation of trading mechanism. While emerging communication and distributed algorithms are also being applied to market bidding [7], it is challenging to implement bidding and power flow calculations in a distributed manner. A bilevel integrated energy sharing mechanism with the decoupled network was designed in [8] for energy hubs and prosumers, in which it did not consider the limitations of information sharing.

To face the dilemma of information sharing limitations and implementation of power flow constraints, this paper proposes a two-level distributed trading mechanism for the regional electricity market, in which distributed energy sources and EH are players, and they are incentivized to participate in the market. At the upper level, the EH is an independent market participant, managing the distributed sources, energy storage, and combined heat and power plant, selling and buying electricity from nodal agents to maximize profit. At the lower level, the nodal agents are divided into supplier and consumer nodes, which as different stakeholders will only pursue the maximization of their interests. The regional system operator (RSO) plays the role of selfless auctioneer, intending to clear the regional electricity market while achieving cross-arbitrage by trading with EH and the upstream grid. The hierarchical electricity trading mechanism is modeled as a mathematical program with equilibrium constraints (MPEC). A distributed iterative algorithm based on the alternating direction method of multipliers (ADMM) is developed to implement distributed clearing process. Case studies based on a realistic regional grid verify the effectiveness of the proposed algorithm and show that the mechanism is effective in improving energy efficiency and reducing the transaction costs.

The remainder of this paper is structured as follows: Section II presents the framework of the two-level market. In Section III, the distributed trading mechanism is proposed. The solution algorithm is proposed in Section IV. The case studies are given in Section V, followed by conclusions in Section VI.

## II. TWO-LEVEL MARKET FRAMEWORK

The flow-based distributed trading market framework is divided into EH-level contains EH, and regional transmission-level which combines consumer, supplier nodal agents and RSO, as shown in Fig. 1. The players in this two-level market consist of EHs, RSO, consumer and supplier nodes. In the regional transmission, the RSO as a coordinator, is trading with EH and the upstream grid while performing DC power flow calculations to clear the regional market (i.e., to find the intersection of offers from suppliers and consumers).

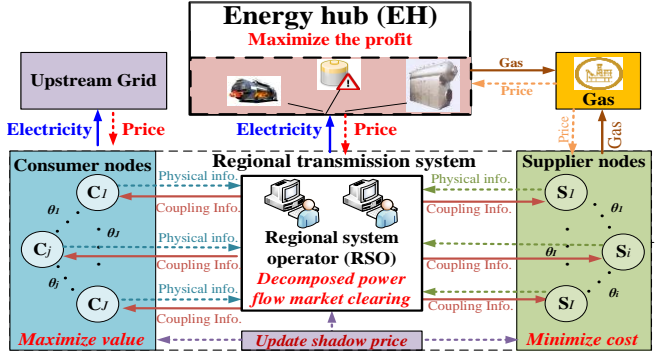

Figure 1. Structure of flow-based distributed trading market.

All nodes $n \in \mathcal{N} = \{1,..,N\}$ in the regional transmission system, are divided into consumer nodes $j \in \mathcal{J} = \{1,..,J\}$ and supplier nodes $i \in \mathcal{I} = \{1,..,I\}$, particularly $\mathcal{J} \cup \mathcal{I} = \mathcal{N}$. Different nodes act as rational nodal agents, competing with other stakeholders to maximize their profits, while sharing as little private physical information as possible, such as cost and value functions. Thus, the proposed flow-based distributed trading mechanism, where players only share voltage angle $\theta$ as a market information, enables clearing the energy market meanwhile implementing DC power flow in a distributed manner. The regional transmission system composed by RSO and all nodal agents as entire stakeholders, cross-arbitrage from the upstream grid and the EH to derive maximum profit. In terms of the EH, it contains energy conversion and storage facilities together, say, combined heat and power (CHP) plant, electrical-powered and electricity storage units. EH determines electricity prices $\lambda_{hb}$ based on trading electricity $P_{hb}$ with regional grid and gas price $\lambda_{gas}$ to maximize the profit.

## III. DECOMPOSED DC FLOW OPTIMIZATION MODEL

The centralized trading model is formulated as an optimization problem presented in (1), which including DC power flow constraints are decomposed by a set of distributed optimization problemes (2)-(5). Each stakeholder solves its dedicated problem which is linked to other problems through some coupling variables and power flow equations as the coupling constraint. The problem of supplier, consumer nodes, and EH are formulated as in (2), (3), and (4) respectively. The dedicated problem for the RSO as the main coordinator of the proposed distributed mechanism is given in (5). Note that proposed distributed power flow calculation each node in the regional grid only needs to share a limited amount of information (voltage angle $\theta_n$). The coupling variables of different nodes are $[P_{n,m}, \theta_n]$, RSO as a selfless coordinator, firstly, aggregates market information, then coordinates, estimates and publishes the market information for the next iteration. It will take several iterations until the convergence condition is met.

### A. Centralized Problem

With a centralized perspective, the regional system trades with the EH and the upstream grid meanwhile performing the DC power flow calculations. The centralized problem is to maximize the profit as given in (1).

$$\min_{\substack{P_{n,g}, P_{n,d}, \\ \varpi_n^s, \varpi_n^b, \\ P_{n,hb}^s, P_{n,hb}^b}} \sum_{n \in N} \begin{pmatrix} c_n(P_{n,g}) - v_n(P_{n,d}) - \lambda_{n,\varpi}^s \varpi_n^s + \lambda_{n,\varpi}^b \varpi_n^b \\ + \lambda_{n,hb}^b P_{n,hb}^b - \lambda_{n,hb}^s P_{n,hb}^s \end{pmatrix} \quad (1a)$$

$$s.t. P_{n,g}^{\min} \leq P_{n,g} \leq P_{n,g}^{\max}, P_{n,d}^{\min} \leq P_{n,d} \leq P_{n,d}^{\max} \quad (1b)$$

$$P_{n,m}^l + \varpi_n^b + P_{n,g} + P_{n,hb}^b = \varpi_n^s + P_{n,d} + P_{n,hb}^s \quad (1c)$$

$$P_{n,m}^l = B_{n,m}(\theta_n - \theta_m), \forall m \in \aleph, m \neq n \quad (1d)$$

$$l_{n,m}^{\min} \leq P_{n,m}^l \leq l_{n,m}^{\max}, \forall m \in \aleph, m \neq n \quad (1e)$$

$$\theta_n^{\min} \leq \theta_n \leq \theta_n^{\max}, \quad (1f)$$

where $c_n$ are the cost of generation, $v_n$ is the value of electricity use, $P_{n,g}, P_{n,d}$ are the generation and demand at node $n$, $\varpi_n^b, \varpi_n^s$ are the electricity purchased from upstream grid and sold to upstream grid, $P_{n,hb}^b, P_{n,hb}^s$ are the electricity purchased from EH and sold to energy hub, $\lambda_{n,\varpi}^s, \lambda_{n,\varpi}^b$ are the wholesale market price of coupling with upstream energy trading $\varpi_n^b, \varpi_n^s$, and $B_{n,m}$ is the susceptance of line $nm$; EH prices $\lambda_{n,hb}^s, \lambda_{n,hb}^b$ are coupled with $P_{n,hb}^b, P_{n,hb}^s$. Here, the objective function contains costs of electricity generation, the value of electricity use, and trading revenue. Constraints (1b) denotes the generation and demand constraints for node $n$. Constraints (1c) - (1f) are the DC power flow model. Constraint (1c) represents the power balance. Constraint (1d) is active power flow through transmission line $l_{n,m}$. Constraint (1e) is the power limits of transmission lines. (1f) indicates the voltage angle constraint.

### B. Supplier Nodal Agents

The objective function of supplier a nodal agent $i \in \mathcal{I} = \{1,...,I\}$ is given in (2a), which consists of the generation cost in the first term, the revenue from electricity sales in the second to the sixth terms and the penalties related to the coupling variable in the seventh term.

$$\min_{\substack{\theta_i, P_{i,g} \\ P_{i,:}^s, P_{i,hb}^s}} \begin{pmatrix} c_i(P_{i,g}) - \lambda_{i,g}(P_{i,g})^T - \lambda_{i,:}^s(P_{i,:}^s)^T - \lambda_{i,hb}^s P_{i,hb}^s \\ -\lambda_{i,\varpi}^s P_{i,\varpi}^s - \lambda_{i,\theta}\theta_i + \frac{\rho}{2}[\|P_{i,:}^s - \hat{P}_{i,:}^s\|_2^2 + \|P_{i,hb}^s - \hat{P}_{i,hb}^s\|_2^2 \\ + \|P_{i,g} - \hat{P}_{i,g}\|_2^2 + \|P_{i,\varpi}^s - \hat{P}_{i,\varpi}^s\|_2^2 + \|\theta_i - \hat{\theta}_i\|_2^2] \end{pmatrix} \quad (2a)$$

$$\text{s.t. } P_{i,g} - P_{i,:}^s \mathbf{1}_{|b|} = 0 \quad (2b)$$

$$P_{i,:}^s = (\theta_i - \hat{\theta})/X_{i,:} \quad (2c)$$

$$P_{i,g}^{\min} \leq P_{i,g} \leq P_{i,g}^{\max}, P_{i,hb}^{s,\min} \leq P_{i,hb}^s \leq P_{i,hb}^{s,\max} \quad (2d)$$

$$\theta_i^{\min} \leq \theta_i \leq \theta_i^{\max}, \quad (2e)$$

where $P_{i,g}$ is the generation at node $i$, $P_{i,:}^s$, $P_{i,hb}^s$ and $P_{i,\varpi}^s$ are the electricity sold by node $i$ to other nodes, energy hub and upstream grid, $\rho$ is the penalty parameter; $\lambda_{i,g}, \lambda_{i,:}^s, \lambda_{i,hb}^s, \lambda_{i,\varpi}^s$ and $\lambda_{i,\theta}$ are the prices of the respective coupling variable $P_{i,g}$, $P_{i,:}^s$, $P_{i,hb}^s$, $P_{i,\varpi}^s$; and $\theta_i$. $\hat{P}_{i,:}^s$, $\hat{\theta}_i$, $\hat{P}_{i,\varpi}^s$ and $\hat{P}_{i,hb}^s$ are the predetermined fixed values of variables, which are determined by RSO and EH. Here, constraint (2b) is the power balance for each supplier node $i$; (2c) is active power flow through transmission line connected with supplier node; (2d) denotes the output constraint for generation; and the constraint of quantity traded with EH; and (2e) indicates the range of voltage angle.

### C. Consumer Nodal Agents

The objective function of a consumer nodal agent $j \in \mathcal{J} = \{1,\ldots,J\}$ is given in (3a):

$$\min_{\substack{P_{j,d},\varpi_j \\ P_{j,hb},P_{:,j}}} \begin{pmatrix} -v_j(P_{j,d}) + \lambda_{j,d}P_{j,d} + (\lambda_{:,j})^T P_{:,j} + \lambda_{j,hb}P_{j,hb} \\ +\lambda_{j,\varpi}\varpi_j + \lambda_{j,\theta}\theta_j + \frac{\rho}{2}[\|P_{:,j} - \hat{P}_{:,j}\|_2^2 + \|P_{j,d} - \hat{P}_{j,d}\|_2^2 \\ + \|P_{j,hb} - \hat{P}_{j,hb}\|_2^2 + \|\varpi_j - \hat{\varpi}_j\|_2^2 + \|\theta_j - \hat{\theta}_j\|_2^2] \end{pmatrix} \quad (3a)$$

$$\text{s.t. } P_{j,d} - \mathbf{1}_{|b|}^T P_{:,j}^b = 0 \quad (3b)$$

$$P_{:,j} = (\hat{\theta} - \theta_j)/X_{:,j} \quad (3c)$$

$$P_{j,d}^{\min} \leq P_{j,d} \leq P_{j,d}^{\max}, P_{j,hb}^{\min} \leq P_{j,hb} \leq P_{j,hb}^{\max} \quad (3d)$$

$$\theta_j^{\min} \leq \theta_j \leq \theta_j^{\max}, \quad (3e)$$

where $P_{j,d}$ is the demand of node $j$; $P_{:,j} = [P_{:,j}^b, P_{:,j}^s]$ is the electricity trading with node $j$ from other nodes; $P_{j,hb} = [P_{j,hb}^b, P_{j,hb}^s]$ is the node $j$'s trading electricity with energy hub; $\varpi_j = [\varpi_j^b, \varpi_j^s]$ is the node $j$'s trading electricity with the upstream grid; $\lambda_{:,j}, \lambda_{j,hb}, \lambda_{j,\varpi}, \lambda_{j,d}$ and $\lambda_{j,\theta}$ are the bid prices of the coupling variables $P_{:,j}$, $P_{j,hb}$, $\varpi_j$, $P_{j,d}$ and $\theta_j$. $\hat{P}_{j,d}, \hat{P}_{:,j}, \hat{\theta}_j, \hat{\varpi}_j$ and $\hat{P}_{j,hb}$ are the predetermined fixed values by RSO and EH. Here, (3b) is the consumer nodes power balance; (3c) is active power flow through transmission line connected with other nodes; (3d) denotes the constraint of demand and quantity traded with EH; (3e) indicates the range of voltage angle.

### D. Energy Hub

Energy hub is a third-party, independent of the RSO, which aims to make a profit by participating in market transactions by solving problem (4).

$$\max_{\substack{\lambda_{hb}^s, \lambda_{hb}^b \\ \hat{P}_{hb}^s, \hat{P}_{hb}^b}} \begin{pmatrix} \lambda_{hb}^s \hat{P}_{hb}^s - \lambda_{hb}^b \hat{P}_{hb}^b - \lambda_{gas} P_{gas} \\ -\frac{\rho}{2}[\|P_{hb}^s - \hat{P}_{hb}^s\|_2^2 + \|P_{hb}^b - \hat{P}_{hb}^b\|_2^2] \end{pmatrix} \quad (4a)$$

$$\text{s.t. } P_{hb,t}^b + P_{gas,t}\eta_e^{chp} + P_{dis,t} - P_{ch,t} = P_{hb,t}^s \quad (4b)$$

$$E_{t+1} = E_t + P_{ch,t}\eta_+^{esu} - P_{dis,t}/\eta_-^{esu} \quad (4c)$$

$$P_{hb,\min}^s \leq P_{hb}^s \leq P_{hb,\max}^s, P_{hb,\min}^b \leq P_{hb}^b \leq P_{hb,\max}^b, \quad (4d)$$

the objective is to maximize the profit of trading, where the first two terms are trading revenue from the transaction with regional grid, the third term represents the cost of purchasing gas, and the last term is the penalty. Constraint (4b) represents the electric balance inside the EH; (4c) is the battery constraint; and the range of EH trading volume is given in (4d).

## IV. DISTRIBUTED SOLUTION ALGORITHMS

In the previous section, a two-level model is devised for the regional grid and EH trading to maximize the profits of all stakeholders. All market players: EH, supplier nodes and consumer nodes within the regional system, as rational stakeholders, aim to maximize revenue by disclosing as little private information as possible in the market transactions. Hence, an ADMM-based algorithm is proposed to clear the market and implement DC power flow, as shown in Fig. 2.

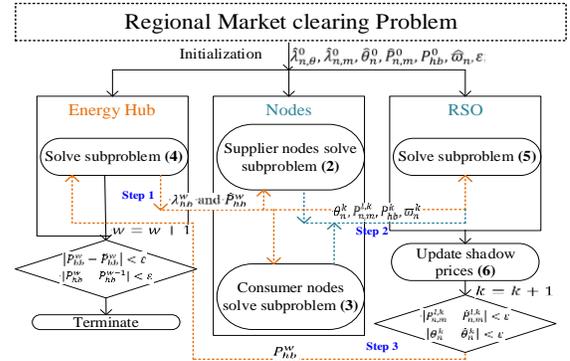

Figure 2. ADMM-based distributed method.

First, EH solves the optimization problem (4) based on initial electricity $P_{hb}^k = [P_{hb}^b, P_{hb}^s]$ to determine the price $\lambda_{hb}^k = [\lambda_{hb}^b, \lambda_{hb}^s]$. All the nodes in the regional system develop bidding strategies, in parallel, consumer nodal agents' strategies $\mathcal{P}_c = [P_{j,d}, P_{:,j}, P_{j,hb}, \varpi_j, \theta_j]$ to maximize value of energy use, and supplier nodal agents pursue minimize cost with strategies $\mathcal{P}_s = [P_{i,g}, P_{i,:}, P_{i,hb}, \varpi_i, \theta_i]$. However, the power flow is coupled in each node and cannot be solved independently, in this case the RSO updates the coupling variables by solving (5) based on the public nodal information.

$$\min_{\substack{\theta_n, \hat{P}_{n,d} \\ \hat{P}_{n,m}, \hat{P}_{n,g}}} \sum_{n \in N} \begin{pmatrix} \lambda_{n,m} \hat{P}_{n,m}^l + \lambda_{n,\theta} \theta - \lambda_{n,g} \hat{P}_{n,g} + \lambda_{n,d} \hat{P}_{n,d} \\ + \frac{\rho}{2} (\|P_{n,m}^l - \hat{P}_{n,m}^l\|_2^2 + \|\theta_n - \hat{\theta}_n\|_2^2 \\ + \|P_{n,g} - \hat{P}_{n,g}\|_2^2 + \|P_{n,d} - \hat{P}_{n,d}\|_2^2) \end{pmatrix} \quad (5a)$$

$$s.t. P_{n,g}^{\min} \le \hat{P}_{n,g} \le P_{n,g}^{\max}, P_{n,d}^{\min} \le \hat{P}_{n,d} \le P_{n,d}^{\max} \quad (5b)$$

$$\hat{P}_{n,m}^l + \hat{\omega}_n^b + \hat{P}_{n,g} + \hat{P}_{n,hb}^b = \hat{\omega}_n^s + \hat{P}_{n,d} + \hat{P}_{n,hb}^s \quad (5c)$$

$$\hat{P}_{n,m}^l = B_{n,m}(\hat{\theta}_n - \hat{\theta}_m), \forall m \in \aleph, m \ne n \quad (5d)$$

$$l_{n,m}^{\min} \le \hat{P}_{n,m}^l \le l_{n,m}^{\max}, \forall m \in \aleph, m \ne n \quad (5e)$$

$$\theta_n^{\min} \le \hat{\theta}_n \le \theta_n^{\max}. \quad (5f)$$

Then, the dual variables are updated according to the latest parameters in (6).

$$\lambda_{n,\theta}^{k+1} = \lambda_{n,\theta}^k + \frac{\rho}{2}[\theta_n^k - \hat{\theta}_n^k], \lambda_{n,m}^{k+1} = \lambda_{n,m}^k + \frac{\rho}{2}[P_{n,m}^k - \hat{P}_{n,m}^k] \quad (6a)$$

$$\lambda_{n,g}^{k+1} = \lambda_{n,g}^k + \frac{\rho}{2}[P_{n,g}^k - \hat{P}_{n,g}^k], \lambda_{n,d}^{k+1} = \lambda_{n,d}^k + \frac{\rho}{2}[P_{n,d}^k - \hat{P}_{n,d}^k] \quad (6b)$$

where $\lambda_{n,\theta} = [\lambda_{i,\theta}, \lambda_{j,\theta}]$, $\lambda_{n,m} = [\lambda_{i,:}, \lambda_{:,j}]$, $\lambda_{n,g} = [\lambda_{i,g}]$ and $\lambda_{n,d} = [\lambda_{j,d}]$ are the shadow prices; $\theta_n = [\theta_i, \theta_j]$ and $P_{n,m} = [P_{i,:}, P_{j,:}]$ are physical variables of consumer and supplier nodes. After several iterations, when the convergence condition is met, algorithm will be terminated. The following assumptions and linearization method are used to find the optimal solution.

## A. Decomposed DC Power Flow

*Assumption 1:* The cost function $c_i: R^i \cup \{+\infty\}$ and value function $v_j: R^j \cup \{+\infty\}$ are closed, proper, and convex.

*Assumption 2:* The overall optimization problem has strictly feasible solutions, and the overall optimization problem constitutes of convex inequality constraints and affine equality constraints.

Assumptions 1 and 2 satisfy the Slater's strong duality condition, so the proposed algorithm has the convergence and global optimality conditions [7]. Therefore, the iterative results of decomposed DC power flow satisfy residual convergence, objective convergence and dual variable convergence.

## B. Binary Expansion Linearisation Method

The product terms $\lambda_{hb}^b \hat{P}_{hb}^b$, $\lambda_{hb}^s \hat{P}_{hb}^s$ in the EH objective function (4a) are nonlinear and non-convex. For these bilinear product terms $\lambda_{hb} \hat{P}_{hb}$, where $\lambda_{hb}$ and $\hat{P}_{hb}$ are continuous variables, a binary expansion method [9] is applied to linearize them. Thus, the binary expansion linearization method is applied to solve EH and transmission system market clearing problem. The objective of EH in (4a) can be reformulated as:

$$\max_{\substack{\lambda_{hb}^s, \lambda_{hb}^b \\ \hat{P}_{hb}^s, \hat{P}_{hb}^b}} \begin{pmatrix} \lambda_{hb}^{s,\min} \hat{P}_{hb}^s + \Delta\lambda_{hb}^s \sum_{\tau=1}^{\Gamma} 2^{\tau-1} \hat{P}_{hb}^s z_\tau^s - \lambda_{gas} P_{gas} - \lambda_{hb}^{b,\min} \hat{P}_{hb}^b \\ -\Delta\lambda_{hb}^b \sum_{\tau=1}^{\Gamma} 2^{\tau-1} \hat{P}_{hb}^b z_\tau^b - \frac{\rho}{2}[\|P_{hb}^s - \hat{P}_{hb}^s\|_2^2 + \|P_{hb}^b - \hat{P}_{hb}^b\|_2^2] \end{pmatrix} \quad (7a)$$

$$s.t. 0 \le \hat{P}_{hb}^s - \hat{P}_{hb}^s z_\tau^s \le \hat{P}_{hb}^{s,\max}(1 - z_\tau^s), \forall \tau \quad (7b)$$

$$0 \le \hat{P}_{hb}^b - \hat{P}_{hb}^b z_\tau^b \le \hat{P}_{hb}^{b,\max}(1 - z_\tau^b), \forall \tau \quad (7c)$$

$$0 \le \hat{P}_{hb}^s z_\tau^s \le \hat{P}_{hb}^{s,\max} z_\tau^s, 0 \le \hat{P}_{hb}^b z_\tau^b \le \hat{P}_{hb}^{b,\max} z_\tau^b, \forall \tau \quad (7d)$$

Cons (4b) - (4e)

where $z_\tau^s, z_\tau^b, \tau = 1, ..., \Gamma$ are binary variables and the step size $\Delta\lambda_{hb}^s$ is $\Delta\lambda_{hb}^s = (\lambda_{hb}^{s,max} - \lambda_{hb}^{s,min})/2^\tau$. The objective (7a) provides a linear expression of $\lambda_{hb}^b \hat{P}_{hb}^b$ and $\lambda_{hb}^s \hat{P}_{hb}^s$. The approximation accuracy of binary expansion can be controlled by the number of expansion segments. The Flow-based distributed trading mechanism in the regional electricity market with EH is a mixed-integer linear programming (MILP) problem, which can be solved by commercial solvers directly.

## V. CASE STUDIES

### A. System Parameters

To verify the effectiveness of the mechanism, we used a case study of regional grid with an EH in Jiangsu Province, in southeastern China. This regional transmission system is connected to the upstream grid at node 1 and node 5, and gas turbines produce electricity at node 4. The EH with CHP and battery storage connects to the regional transmission system at node 2, as shown in Fig. 3. Detailed system data can be found in [10]. The price of gas is 26$/MWh and the price of selling electricity to the upstream grid is 90$/MWh.

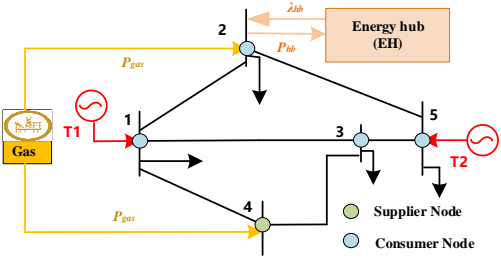

Figure 3. The structure of the test node system.

### B. Analysis of Regional Electricity Market Clearing

At the lower level, the regional system participates in the wholesale market by interacting with the upstream grid while trading with EH. Fig. 4 shows the energy trading results and Fig. 5 shows the price in the market. In the day-ahead, nodal agents in regional transmission system buy electricity from the upstream grid via node 1 while reselling electricity to the upstream grid via node 5 to achieve the cross-arbitrage, due to the difference between the purchase and the sale. During periods 7-8, the price of electricity sold by EH is in the valley periods and is almost equal to the price of electricity sold from the upstream grid. Hence, the transmission system buys part of the electricity from EH and resells it to the upstream grid through node 5. In addition, the prices of electricity sold by EH and by the upstream grid are very close in periods 19 to 22, while prices are not the lowest at this point, the regional system can still cross arbitrage in transactions between EH and the upper grid.

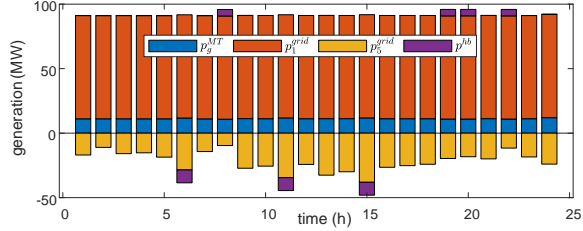

Figure 4. The results of the energy trading of the transmission system.

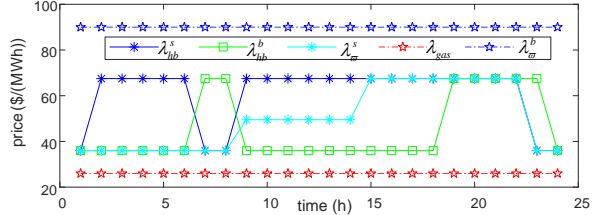

Figure 5. The prices in the two-layer trading market.

The amounts of purchased gas by EH and state of charge of the battery storage inside the EH are indicated in Fig.6. The EH purchases gas before discharging to benefit from selling electricity to the regional system. It is interesting to note that during periods 10-15, EH only buys electricity from regional transmission system instead of purchasing gas and converting it into electricity. It is because the cost of electricity purchased by EH from the regional grid during this period is much lower than the cost of purchasing gas for gas-to-electric conversion.

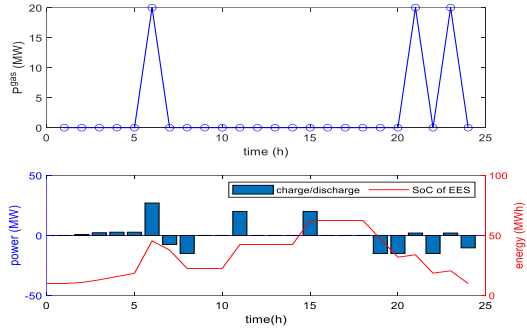

Figure 6. The amounts of purchased gas and state of charge of the battery storage in the EH.

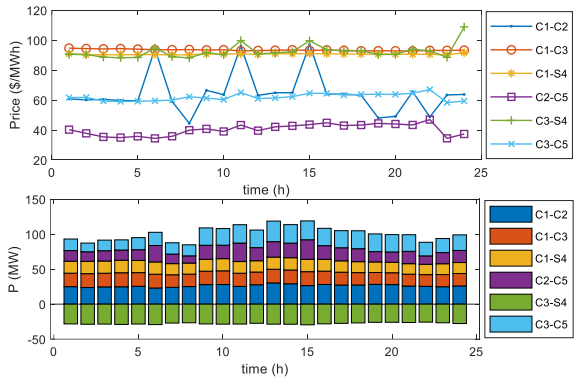

Figure 7. Price and electricity of transactions between different nodes.

It can be observed from Fig. 7 that nodes C1 and C2 have the most variation in price, while C2 has the option to purchase power from C1 and EH, both of which fluctuate at around 60$/MWh. In addition, the price of C1-C2 varies proportionally with the quantity of transactions. The prices of C1-C3, C1-S4 and C3-S4 are almost the highest, because C1 purchases energy from the upstream grid and resells it to C3 and S4 to increase the revenue. It is interesting to note that these three prices are very close to each other. C3 as a consumer node, who can only obtain the energy from S4 and C1, meanwhile every node wants to benefit from the market that leads to this price phenomenon. Nodes C2-C5 have the lowest prices, and C2 can sell electricity to EH and node C5, simultaneously. Comparing the EH purchase price in Fig. 5 with the purchase price for C5 in Fig. 7, the two prices are comparable in periods 0 to 5 and 9 to 18. However, the purchase prices of EH increase dramatically during the periods 18 to 23, meanwhile the prices of C2-C5 have also risen, steadily. Thus, during this period, C2 sells electricity to C5 only, not to the EH. Next, the results of the proposed distributed mechanism is compared with the centralized method. Table I shows the comparison of players' value (v), cost (c), revenue and social welfare (SC) in the centralized and distributed method. EH and regional transmission systems (RS) have more benefits in the distributed clearing approach than in the centralized clearing method.

TABLE I COMPARISON OF RESULTS

| Model | Player | υ/c($) | Revenue ($) | SC ($) |
|---|---|---|---|---|
| Distributed | C1 | 408.09 | 1344.15 | 363.06 |
|  | C2 | 650.67 | 431.63 |  |
|  | C3 | 373.34 | -376.23 |  |
|  | S4 | 1123.87 | -837.67 | 994.96 |
|  | C5 | 468.94 | 433.07 |  |
|  | EH | 414.11 | 12.5 |  |
| Centralized | C1 | 213.42 | 19.47 | 163.71 |
|  | C2 | 703.26 | 702.94 |  |
|  | C3 | 111.05 | 111.05 |  |
|  | S4 | 1201.77 | -1201.77 | 352.26 |
|  | C5 | 468.95 | 720.57 |  |
|  | EH | 131.2 | 3.2 |  |

This is because under the distributed algorithm, nodes C1 and C5 act as rational players and increase their revenue by trading with the upstream grid and reselling to other nodes.

### C. Performance of algorithms

The convergence process of the flow-based distributed algorithm is illustrated in Fig. 8.

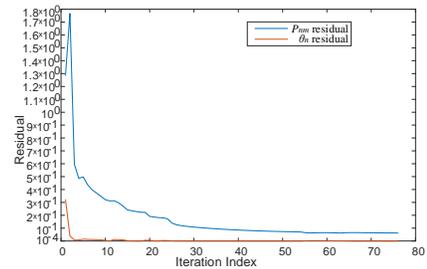

Figure 8. The convergence process of the flow-based distributed algorithm.

The residual variables, as shown in Fig. 2, are $P_{nm}$ and $\theta_n$. The stopping criteria are $\epsilon^p = 10^{-2}$ and $\epsilon^\theta = 10^{-3}$. Each iteration takes 0.69s, and the total computation time is 52s, satisfying the conditions for day-ahead market clearing. When the number of iterations is less than ten, there is tremendous shock, but when the number of iterations reaches 30, the convergence curve gradually stabilizes.

## VI. Conclusion

This paper proposes a flow-based hierarchical trading mechanism that combines the regional transmission system and EH within a regional electricity market. At the upper level, EH participates as a stakeholder purchasing fuel from the gas market and trading with the nodal agents in regional transmission system; whereas at the lower level, different nodal agents intend to cross arbitrage from EH and upstream grid. The interaction between nodal agents and EH is formulated through an MPEC, for which an ADMM-based distributed algorithm is proposed to realize the regional electricity market clearing. The DC power flow is calculated by the RSO and nodes in a distributed manner while considering the privacy concerns. Case studies show that the flow-based trading mechanism provides a good incentive for EH and distributed energy resources to participate in the market to improve energy efficiency. Furthermore, the effectiveness and convergence of the algorithm are also well proved.

## References

[1] J. Wang et al., "Exploring the trade-offs between electric heating policy and carbon mitigation in China," *Nat. Commun.*, vol. 11, no. 1, pp. 6054, 2020.
[2] Y. Li, Z. Li, F. Wen, M, Shahidehpour, "Privacy-preserving optimal dispatch for an integrated power distribution and natural gas system in networked energy hubs", *IEEE Trans. Sustain. Energy*, vol. 10, no. 4, pp. 2028-2038, 2019.
[3] P. Li, W. Sheng, Q. Duan, Z. Li, C. Zhu and X. Zhang, "A Lyapunov optimization-based energy management strategy for energy hub with energy router," *IEEE Trans. Smart Grid*, vol. 11, no. 6, pp. 4860-4870, 2020.
[4] D. Wang, Q. Hu, H. Jia, K. Hou, W. Du, N. Chen, X. Wang, M. Fan, "Integrated demand response in district electricity-heating network considering double auction retail energy market based on demand-side energy stations", *Appl. Energy*, vol. 248, pp. 656-678, 2019.
[5] V. Davatgaran, M. Saniei, and S. S. Mortazavi, "Optimal bidding strategy for an energy hub in energy market, " *Energy*, vol. 148, pp. 482–493, 2018.
[6] L. P. M. I. Sampath, A. Paudel, H. D. Nguyen, E. Y. S. Foo and H. B. Gooi, "Peer-to-peer energy trading enabled optimal decentralized operation of smart distribution grids," *IEEE Trans. Smart Grid*, vol. 13, no. 1, pp. 654-666, 2022
[7] S. Boyd, N. Parikh, E. Chu, B. Peleato, and J. Eckstein, "Distributed optimization and statistical learning via the alternating direction method of multipliers," *Found. Trends Mach. Learn.*, vol. 3, no. 1, pp. 1–122, 2011.
[8] N. Liu, L. Tan, H. Sun, Z. Zhou and B. Guo, "Bilevel heat–electricity energy sharing for integrated energy systems with energy hubs and prosumers," *IEEE Trans. Ind. Informat.*, vol. 18, no. 6, pp. 3754-3765, June 2022
[9] M. V. Pereira, S. Granville, M. H. Fampa, R. Dix, and L. A. Barroso, "Strategic bidding under uncertainty: a binary expansion approach," *IEEE Trans. Power Syst.*, vol. 20, no. 1, pp. 180–188, 2005.
[10] L. Wang, "A flow-based distributed trading mechanism in regional electricity market with energy hub", 2022, https://www.dropbox.com/s/d1qojn5kulxigki/EEM22_full%20vision%20with%20param.docx?dl=0.